# A New Differential Test for Series of Positive Terms

*Yi-Fang Chang*
*Department of Physics, Yunnan University, Kunming 650091, China*
(e-mail: yifangchang1030@hotmail.com)

**ABSTRACT** A new differential test for series of positive terms is proved. Let *f(x)* be a positive continuous function corresponded to a series of positive terms $\sum_{k=1}^{\infty} f(k)$, and *g(x)* is a derivative of reciprocal of *f(x)*, i.e., $\frac{d}{dx}[\frac{1}{f(x)}] = g(x)$. Then, if $fgx \geq 1+\alpha(\alpha>1)$ for enough large x, the series converges; if $fgx \leq 1$, the series diverges. The rest may make the limit form, and is universal and complete.

The convergence or divergence for a series of positive terms is elementary in calculus, and has been many tests. Besides the basic comparison test, there are the D'Alembert ratio rest, the integral test [1], the Cauchy root test, the Raabe test, the Kummer test, the Abel-Dini test, etc [2,3]. But, the applicable regions of these tests are not usually the same, and are generally finite. Using a similar method with the integral test, we propose a new differential test for series of positive terms, and discuss some examples.

The differential test. Let $\sum_{k=1}^{\infty} f(k)$ be a series of positive terms, *f(x)* is a corresponding positive continuous function, and *g(x)* is a derivative of reciprocal of *f(x)*, i.e., $\frac{d}{dx}[\frac{1}{f(x)}] = g(x)$. Then, if $fgx \geq 1+\alpha(\alpha>0)$ for enough large x, the series converges; if $fgx \leq 1$ the series diverges.

Proof Since $g(x) = \frac{d}{dx}[\frac{1}{f(x)}] = -\frac{d}{dx}[\ln f(x)]\frac{1}{f(x)}$, if $fg = -\frac{d}{dx}[\ln f(x)] \geq \frac{1+\alpha}{x}$, then the integral from 1 to x is obtained $-\ln\frac{f(x)}{f(1)} \geq (1+\alpha)\ln x$, $\therefore f(x) \leq \frac{f(1)}{x^{1+\alpha}}$, the series $\left\{\frac{1}{k^{1+\alpha}}\right\}$ and {*f(k)*} converge; if $fg = -\frac{d}{dx}[\ln f(x)] \leq \frac{1}{x}$, the integral is obtained $\ln\frac{f(1)}{f(x)} \le \ln x$, $\therefore\ f(x) \ge \frac{f(1)}{x}$, the series $\left\{\frac{1}{k}\right\}$ and {*f(k)*} diverge.

In calculus the differentiation is simple, and can apply to various composite functions of any elementary functions, therefore, the test is also very simple, and the applicable region is wide. Moreover, a calculating result for any series must be $fgx \ge 1+\alpha(\alpha > 0)$ or $fgx \le 1$, so the rest should be universal and complete.

If $f$ is a discrete function $a_n$, the difference is substituted for differentiation

$$g = \frac{a_{n+1}^{-1} - a_n^{-1}}{1} = \frac{a_n - a_{n+1}}{a_n a_{n+1}},\quad \lim_{n\to\infty} nfg = \lim_{n\to\infty} n(\frac{a_n}{a_{n+1}} - 1) = c$$, i.e., the Raabe test.

In many cases, the test may make the limit form, i.e., *fgx*→*c* as $x \to \infty$. Then the series converges if *c*>1 and diverges if *c*<1. Therefore, it may combine the L’¯ Hbpit a’¡¯s rul

Example 1. The series $\sum_{n=1}^{\infty} f(n)$ and $f(n) = \frac{1}{n^p}$, then $f(x) = \frac{1}{x^p}$,

$g(x) = px^{p-1}$, $fgx = p, \therefore$ the series converges if *p*>1, and diverges if $p \le 1$.

Example 2. The series $f(n) = a^n b^{n-1} + a^n b^n$ cannot apply the ratio test. Let

$f(x) = a^x b^{x-1} + a^x b^x$, $g(x) = -\frac{\ln(ab)}{a^x b^{x-1}(1+b)}$, $fgx = -x\ln(ab), \therefore$ the series converges if *ab*<1, and diverges if $ab \ge 1$.

Example 3. The series $f(n) = \frac{1}{n!}(na)^n, (a > 0)$ must apply the Raabe test. Let

$f(x) = \frac{1}{\Gamma(x+1)}(xa)^x$, $g(x) = (xa)^{-x}[\Gamma'(x+1) - (1+\ln ax)\Gamma(x+1)]$. When x is very large number, $\Gamma'(x+1) \approx [(\frac{x}{e})^x \sqrt{2\pi x}]' = \Gamma(x+1)(\frac{1}{2x} + \ln x)$, $fgx = \frac{1}{2} - x\ln(ae), \therefore$

$\lim_{x\to\infty} fgx = \infty$ for *a*<1/e, the series converges; fgx$\le$1/2 for $a \ge$1/e the series diverges.

Example 4. The series $f(n) = \frac{1}{n^p}\sin\frac{\pi}{n}$, let $f(x) = \frac{1}{x^p}\sin\frac{\pi}{x}$,

$g(x) = \frac{x^{p-1}}{\sin(\pi/x)}[p + \frac{\pi}{x}\mathrm{ctg}\frac{\pi}{\mathrm{x}}]$, $fgx = p + \frac{\pi}{x}ctg\frac{\pi}{x}$, $\lim_{x\to\infty}\frac{\pi}{x}ctg\frac{\pi}{x} = \lim_{x\to\infty}\cos^2\frac{\pi}{x}$ =1,

fgx=1+*p*, the series converges for *p*>0, and diverges for $p \le 0$.

Further, the differentiation test may be applied to test for the general series of functional terms, and test for the infinite integral whose limits of integral from positive to infinite. The convergence or divergence of the integral is determined by the differential of the integrand.

Example 5. The integral $\int_1^{\infty} \frac{dx}{x^p + x^q}$ , $f(x) = \frac{1}{x^p + x^q}$ , $g(x) = px^{p-1} + qx^{q-1}$ ,

$fgx = \frac{px^p + qx^q}{x^p + x^q}$. Using the L’¯ Hopita’¡¯s ru, fgx$\rightarrow$c=$p$ for $p>q$ or c=$q$ for $p<q$, so max ($p$,$q$)>1, the infinite integral converges.

## References

1.L.H.Loomis, Calculus (Third Edition). Addison-Wesley Publishing Company. 1982.

2.G.H.Hardy, A Course of Pure Mathematics, Tenth Edition. Cambridge Mathematical Library. 1993.

3.T.J.I’A.Bromwich, An Introduction to the Theory of Infinite Series. AMS Chelsea Publishing. 2005.